\documentclass[11pt]{amsart}
\usepackage{amssymb,amsmath,amsthm}  %

\chardef\bslash=`\\ 

\newtheorem*{thm}{Theorem} 

\newtheorem{prop}{Proposition}

\theoremstyle{definition}
\newtheorem{defn}{Definition}[]

\theoremstyle{remark}
\newtheorem*{rem}{Remark} 
\newtheorem*{notation}{Notation}

\newcommand{\propref}[1]{Proposition~\ref{#1}}

\newcommand{\Na}{\mathbb{N}}
\newcommand{\fC}{\mathfrak{C}}
\newcommand{\fCz}{{\mathfrak{C}}_0}
\newcommand{\fCl}{{\mathfrak{C}}_\lambda}

\newcommand{\Zi}{\mathbb{Z}}

\newcommand{\R}{\mathbb{R}}
\newcommand{\fR}{\mathfrak{R}}

\newcommand{\C}{\mathbb{C}}
\let\bsy\boldsymbol
\def\P{\bsy P}
\def\T{\bsy T}
\def\t{\bsy t}

\def\p{\bsy p}
\def\f{\bsy f}
\def\K{\bsy K}
\def\k{\bsy k}

\def\F{\bsy F}

\def\H{\mathcal{H}}
\newcommand{\A}{\mathcal{A}}

\newcommand{\Ss}{\mathcal{S}}

\DeclareMathOperator{\sgn}{sgn}

\newcommand{\eval}[2][\right]{\relax
  \ifx#1\right\relax \left.\fi#2#1\rvert}

\let\abs=\envert

\let\norm=\enVert

\let\hnorm=\henVert
\newcommand{\kskob}[1]{\left[#1\right]}
\let\ksk=\kskob
\newcommand{\skob}[1]{\left(#1\right)}
\let\sk=\skob
\newcommand{\skoba}[1]{\left\langle#1\right\rangle}
\let\ska=\skoba

\let\hska=\hskoba

\begin{document}
\renewcommand{\sectionmark}[1]{}
\title[Integral representations]{Integral representations \\ of unbounded operators \\ by arbitrarily smooth
Carleman kernels}
\author[I. M. Novitski\u i]{Igor M. Novitski\u i$^*$}
\curraddr{Institute for Applied Mathematics, Russian Academy of Sciences,
9, Shevchenko Street, Khabarovsk 680 000, Russia}
\email{novim@iam.khv.ru}
\dedicatory{Dedicated to Professor V.\ B.\ Korotkov in celebration
               of the 65th anniversary of his birth}
\thanks{$^*$Research supported in part by Deutscher Akademischer Austauschdienst
grants A/00/06533, A/01/06997. A version \cite{Je} of this paper was written
when the author enjoyed the hospitality of the
Mathematical Institute of Friedrich-Schiller-University, Jena, Germany}
\keywords{Integral operator, Carleman operator, Hilbert-Schmidt operator,
limit spectrum, unitary equivalence}
\subjclass{47G10, 45P05}
\begin{abstract}
In this paper, we give a characterization of all closed linear operators
in a separable Hilbert space which are unitarily equivalent to an integral
operator in $L_2(\R)$ with bounded and arbitrarily smooth Carleman kernel
on $\R^2$. In addition, we give an explicit
construction of corresponding unitary operators.
\end{abstract}
\maketitle

\section{Introduction}
Throughout, $\H$  will denote a separable Hilbert space,
$\fC(\H)$ the set of all closed linear operators densely defined in $\H$,
$\fR(\H)$ the algebra of all bounded linear
operators on $\H$, and $\C$, and $\Na$, and $\Zi$, the complex
plane, the set of all positive integers, the set of all integers,
respectively.
For an operator $S$ in $\fC(\H)$, $S^*$ will denote the Hilbert space
adjoint of $S$.

Let $S:D_S\to\H$ be an operator of $\fC(\H)$.
A number $\lambda\in\C$ is said to
belong to the {\it limit spectrum\/} $\sigma_c(S)$ of $S$
if there exists an orthonormal sequence $\{e_n\}\subset D_S$ such that
$\lim\limits_{n\to\infty}\hnorm{(S-\lambda)e_n}=0$.
Let $\fCl(\H)$ denote the collection of all operators $S\in\fC(\H)$
for which $\lambda\in\sigma_c(S^*)$.

Let $\R$ be the real line $(-\infty,+\infty)$ equipped with the Lebesgue measure,
and let $L_2=L_2(\R)$ be the Hilbert
space of (equivalence classes of) measurable complex--valued functions on
$\R$ equipped with the inner product
$$
\ska{f,g}=\int_{\R} f(s)\overline{g(s)}\,ds
$$
and the norm
$\Vert f\Vert=\ska{f,f}^{\frac{1}2}$.
A linear operator
                 $T : D_T \to  L_2$,
                 where the domain
                       $D_T$
                       is a dense linear manifold in
                                                    $L_2$,
is said to be {\it integral\/}
if there exists a measurable function
$\T$ on $\R^2$, a {\it kernel\/},
such that, for every $f\in D_T$,
$$
               (Tf)(s)=\int_{\R} \T(s,t)f(t)\,dt
$$
for almost every $s$ in $\R$.
If $T$ is an integral operator  such that $D_T=L_2$ then
$T\in\fR(L_2)$ (see \cite[Theorem I.2.8]{Kor:book1}, \cite[Theorem 3.10]{Halmos:Sun}).
A kernel $\T$ on
$\R^2$ is said to be {\it Carleman\/} if $\T(s,\cdot) \in L_2$
for almost every fixed $s$ in $\R$.
An integral operator with a kernel $\T$ is
called {\it Carleman\/} if $\T$ is a Carleman kernel, and it is called
{\it bi-Carleman\/} if both $\T$ and $\T^*$ ($\T^*(s,t)=\overline{\T(t,s)}$)
are Carleman kernels.
Every Carleman kernel, $\T$, defines a {\it Carleman
function\/} $\t$ from $\R$ to $L_2$ by
$\t(s)=\overline{\T(s,\cdot)}$
for all $s$ in ${\R}$ for which $\T(s,\cdot)\in L_2$.

We shall recall a characterization of Carleman representable operators \cite[Theorem IV.3.6]{Kor:book1}:
\begin{prop}\label{Kor}
A necessary and sufficient condition that
an operator $S\in\fC(\H)$ be unitarily equivalent to an integral operator with
Carleman kernel is
that $S$ belong to $\fCz(\H)$.
\end{prop}
While this result deals with Carleman operators in general,
and the question whether their kernels can have certain analytic properties
is largely beside the point, given any non-negative integer $m$,
we impose on a Carleman kernel $\K$ the following smoothness conditions:
\begin{enumerate}
\renewcommand{\labelenumi}{(\roman{enumi})}
\item the function $\K$ and all its partial derivatives
on $\R^2$ up to order $m$ are in $C(\R^2,\C)$,
\item the Carleman function $\k$,
$\k(s)=\overline{\K(s,\cdot)}$,
and all its (strong) derivatives on ${\R}$ up to order $m$
are in $C(\R,L_2)$.
\end{enumerate}
Here and throughout  $C(X,B)$, where $B$ is a Banach space (with norm
$\Vert \cdot\Vert _B$), denote the Banach space (with the norm $\Vert f\Vert _{C(X,B)}
=\sup\limits_{x\in X}\,\Vert f(x)\Vert _B$) of
continuous $B$-valued functions defined on a
locally compact space $X$ and {\it vanishing at infinity\/} (that is, given any
$f \in C(X,B)$ and $\varepsilon>0$, there exists a compact subset
$X(\varepsilon,f) \subset X$ such that $\Vert f(x)\Vert _{B} < \varepsilon$
whenever $x \not\in  X(\varepsilon,f)$).

\begin{defn} A function $\K$ that satisfies Conditions (i), (ii)
is called a {\it $SK^m$-kernel\/} \cite{nov:Isra}. In addition, a $SK^m$-kernel $\K$ is
called a {\it $K^m$-kernel\/} (see \cite{Nov:Lon}, \cite{nov:91}) if
the conjugate transpose function $\K^*$,
$\K^*(s,t)=\overline{\K(t,s)}$,
is also a $SK^m$-kernel, that is,
its Carleman function $\k^*$, $\k^*(s)
=\overline{\K^*(s,\cdot)}$, satisfies Condition (ii).
\end{defn}
We are now in a position to formulate a series of known results on
integral representations of operators by $SK^m$-kernels.
\begin{prop}\label{msmooth}
Let $m$ be a fixed non-negative integer. Then
\begin{enumerate}
\renewcommand{\labelenumi}{\bf(\Alph{enumi})}
\item
if for an operator $S\in\fCz(\H)$
there exist a dense in $\H$ linear manifold $D$
and an orthonormal sequence $\{e_n\}$
such that
\begin{equation*}\label{biCarleman}
\{e_n\}\subset D\subset D_S\cap D_{S^*},\quad
   \lim\limits_{n\to\infty}\hnorm{Se_n}=0,\quad
   \lim\limits_{n\to\infty}\hnorm{S^* e_n}=0,
\end{equation*}
then there exists
an explicitly computable  unitary operator
$U_m:\H\to L_2$ such that $T=U_mSU_m^{-1}$ is a bi-Carleman operator having a
$K^m$-kernel {\rm(see \cite{nov:91}, \cite{nov:92})},
\item
if for an operator family $\left\{S_\alpha:\alpha \in \A\right\}\subset
\fCz(\H)\cap\fR(\H)$
there exists an orthonormal sequence
$\left\{e_n\right\}$
such that
$$
\lim\limits_{n\to\infty}\sup\limits_{\alpha\in\A}\hnorm{S_\alpha^*e_n}=0,
\quad
\lim_{n\to\infty}\sup_{\alpha\in\A}\hnorm{S_\alpha e_n}=0,
$$
then there exists an explicitly computable unitary operator $U_m:\H\to L_2$ such that all the operators
$U_mS_\alpha U_m^{-1}$ $(\alpha\in\A)$
are bounded bi-Carleman operators having $K^m$-kernels {\rm(see \cite{Nov:Lon})},

\item
if for a countable family $\left\{S_r:r\in\Na\right\}\subset\fCz(\H)\cap\fR(\H)$
there exists an orthonormal sequence
$\left\{e_n\right\}$
such that
$$
\sup_{r\in\Na}\hnorm{S^*_r e_{n}}\rightarrow0\quad
 \text{as\ } n \rightarrow \infty,
$$
then there exists
an explicitly computable unitary operator $U_m:\H\to L_2$ such that all the operators
$U_mS_r U_m^{-1}$ $(r\in \Na)$ are bounded Carleman operators having
$SK^m$-kernels \rm{(see \cite{nov:Isra})}.
\end{enumerate}
\end{prop}

In the proposition just formulated, the ``explicit computability''
of the unitary operators $U_m$ means that they are given by
\begin{equation}\label{Um}
U_mf_n=u_n\quad (n\in\Na),
\end{equation}
where $\{f_n\}$, $\{u_n\}$
are orthonormal bases in $\H$ and $L_2$, respectively.
As the basis $\{f_n\}$, the proof of Proposition uses an arbitrary
orthonormal basis in $\H$ provided that $\{e_n\}\subset\{f_n\}$.
Given a non-negative integer $m$,
the general method to obtain the second desired basis $\{u_n\}$
for $L_2$ is to rearrange in a special way the functions
\begin{equation}\label{mbasis}
u_{k,l}(s)=\begin{cases}
\dfrac1{\sqrt{|t_{l+1}-t_l|}} y_{k}\skob{\dfrac{s-t_l}{|t_{l+1}-t_l|}}
&\text{if $s\in\kskob{t_l,t_{l+1}}$},\\
0&\text{if $s\in\R\setminus\kskob{t_l,t_{l+1}}$}
\end{cases}\quad(k\in\Na,\ l\in\Zi),
\end{equation}
where $-\infty<\dots<t_l<t_{l+1}<\dots<+\infty$, with $t_l\to\pm\infty$ as
$l\to\pm\infty$, and
$\{y_k\}$ is an orthonormal basis in $L_2[0,1]$ whose terms are
eigenfunctions  of a self-adjoint differential operator $L$ on $[0,1]$
induced by the differential expression $l(y)=y^{2m+2}$ and by the boundary
conditions $y^{(i)}(0)=y^{(i)}(1)=0$ ($i=0,\dots,m$).

The purpose of this paper is to restrict the conclusion of \propref{Kor}
to {\it arbitrarily\/} smooth Carleman kernels. Now define these kernels.

\begin{defn}
We say  that a function $\K$ is a {\it $SK^\infty$-kernel}
({\it $K^\infty$-kernel}) if
it is a $SK^m$-kernel ($K^m$-kernel) for each non-negative integer $m$.
\end{defn}

\begin{thm}\label{thm-main}
If $S\in\fCz(\H)$,
then there exists a unitary operator $U_\infty:\H\to L_2$ such that the
operator $T=U_\infty SU^{-1}_\infty$ is a  Carleman operator having $SK^\infty$-kernel.
\end{thm}

It should be noted that this theorem, which is our main result, both extends
the singleton case of the statement {\bf(C)} of \propref{msmooth}
to unbounded operators of $\fCz(\H)$ and restricts its conclusion
to $SK^\infty$-kernels.

Theorem will be proved in the next section of the present paper.
The proof yields an explicit construction of the unitary operator
$U_\infty:\H\to L_2$. The construction of $U_\infty$ is
independent of those spectral points of $S$ that are different from $0$,
and is defined, as in \eqref{Um}, by $U_\infty f_n=u_n$ ($n\in\Na$),
where $\{f_n\}$, $\{u_n\}$ are orthonormal bases in $\H$ and $L_2$, respectively,
whose elements can be explicitly described in terms of the operator
$S$. Unlike the functions \eqref{mbasis}, the functions $u_n$ ($n\in\Na$)
in the proof below are naturally required to be arbitrarily smooth on $\R$,
and Section~3 of this paper exhibits an explicit example of such a basis $\{u_n\}$ adopted
from the wavelet theory.

         \section{Proof of Theorem}
The proof consists of three steps. The first step is some preparation for the next
two steps. In this step we split the operator $S\in\fC_0(\H)$ in order to construct
auxiliary operators $J$, $\varGamma$, $B$, $Q$. Using these operators,
in the second step we describe suitable orthonormal bases $\{u_n\}$, $\{f_n\}$
in $L_2$, $\H$, respectively, and suggest (as $U_\infty$ in the theorem)
a certain unitary operator from $\H$ to $L_2$ that sends the basis
$\{f_n\}$ onto the basis $\{u_n\}$. The rest of the proof is a straightforward verification
that the constructed unitary operator is indeed as desired.

\subsection*{Preparing} Let $S\in\fCz(\H)$, and let $\{e_k\}\subset D_{S^*}$
be that orthonormal sequence for which
\begin{equation}\label{null-seq}
\sum_k\hnorm{S^*e_k}^{\frac14}<\infty
\end{equation}
will do (the sum notation $\sum\limits_k$ will always be used instead of
the more detailed symbol $\sum\limits_{k=1}^\infty$).
Let $H$ be the closed linear span of the $e_k$'s, and let
$H^\perp$ be the orthogonal complement of $H$ in $\H$.
Since $S^*\in\fC(\H)$, we have
\begin{equation}\label{hvd}
H\subset D_{S^*}.
\end{equation}
If $E$ is the orthogonal projection onto $H$, consider the decomposition
\begin{equation}\label{splitting}
  S=(1-E)S+ES.
\end{equation}
Since $E\in\fR(\H)$, it follows via definition of the adjoint that
$(ES)^*=S^*E$.
Observe that the operator $J=S^*E$ is nuclear; this property
follows immediately from \eqref{null-seq} and \eqref{hvd}.
 Assume, with no loss of generality, that
$\dim H^\perp=\dim H=\infty$, and choose an orthonormal basis
$\left\{e_k^\perp\right\}$ in $H^\perp$ so that
\begin{equation}\label{hkvd}
  \left\{e_k^\perp\right\}\subset (1-E)D_{S^*}. 
\end{equation}
For each $f\in D_{S^*}$ and for each $h\in \H$, let
\begin{equation}\label{defdknums}
z(f)=\hnorm{S^*f}+1, \quad d(h)=\hnorm{Jh}^\frac{1}4+\hnorm{J^*h}^\frac{1}4+
\hnorm{\varGamma^*h},
\end{equation}
where
\begin{equation}\label{lamgam}
\varGamma=S^*\Lambda, \text{\ and\ } \Lambda =\sum_k\frac 1{kz\sk{e_k^\perp}}
\hska{\cdot,e_k^\perp} e_k^\perp.
\end{equation}
Since, for the basis $\{f_n\}=\{e_k\}\cup\left\{e_k^\perp\right\}$,
$$
\sum_n\hnorm{S^*\Lambda f_n}^2= \sum_k\dfrac{\hnorm{S^*e_k^\perp}^2}
{k^2(\hnorm{S^*e_k^\perp}+1)^2}\le\frac{\pi^2}6,
$$
it follows that $\varGamma^*$ is a Hilbert-Schmidt operator.

If $J=\sum_n s_n\hska{\cdot,p_n} q_n$ is the Schmidt decomposition for $J$,
then the closedness of $S$ implies that, for every $g\in D_S$,
\begin{equation}\label{ES}
  ESg=(ES)^{**}g=(S^*E)^*g=J^*g=\sum\limits_n s_n\hska{g,q_n} p_n;
\end{equation}
here the $s_n$ are the singular values of $J$ (eigenvalues of
$\left(JJ^*\right)^{\frac1{2}}$), $\left\{p_n\right\}$, $\left\{q_n\right\}$  are
orthonormal sets (the $p_n$ are eigenfunctions for $J^*J$ and $q_n$ are
eigenfunctions for $JJ^*$).

Introduce one more auxiliary operator $B$ by
\begin{equation}\label{B}
B=\sum_n s_{n}^{\frac1{4}}\hska{\cdot,p_n} q_{n},
\end{equation}
and observe that, by the Schwarz inequality,
\begin{equation}\label{schw}
\hnorm{Bf}=\hnorm{\left(J^*J\right)^{\frac18}f}\leq\hnorm{J f}^{\frac14},\quad
\hnorm{B^*f}=\hnorm{\left(JJ^*\right)^{\frac18}f}\leq \hnorm{J^*f}^{\frac14}
\end{equation}
if $\norm{f}=1$.

Now consider the operator $Q=(1-E)S$. The property \eqref{hkvd}
provides the representation
\begin{equation}\label{qoper}
Qg=\sum\limits_k\hska{Qg,e_k^\perp} e_k^\perp=
  \sum\limits_k\hska{g,S^*e_k^\perp} e_k^\perp\quad\text{for all $g\in D_S$}.
\end{equation}

\subsection*{Defining a unitary $\bsy {U_\infty}$} This step is to  construct a candidate for the desired
unitary operator $U_\infty$ in the theorem.

\begin{notation}
If an equivalence class
$f\in L_2$ contains a function belonging to $C(\R,\C)$, then we shall use
$\ksk{f}$ to denote that function.
\end{notation}

Take any orthonormal basis $\{u_n\}$ for $L_2$ which satisfies conditions:
\begin{enumerate}
\renewcommand{\labelenumi}{(\alph{enumi})}
\item the terms of the derivative sequence
$\left\{\ksk{u_n} ^{(i)}\right\}$ are in $C(\R,\C)$, for each $i$
(here and throughout, the letter $i$ is reserved for integers in $[0,+\infty)$),

\item $\{u_n\}=\{g_k\}_{k=1}^\infty\cup\{h_k\}_{k=1}^\infty$, where
$\{g_k\}_{k=1}^\infty\cap\{h_k\}_{k=1}^\infty=\varnothing$,
and,  for each $i$,
\begin{equation}\label{hki}
\sum_k H_{k,i}<\infty\quad
\text{with $H_{k,i}=\norm{\ksk{h_k}^{(i)}}_{C(\R,\C)}$}\quad (k\in\Na),
\end{equation}
\item there exist a subsequence
$\left\{x_k\right\}_{k=1}^\infty\subset\{e_k\}$ and
a strictly increasing sequence $\left\{n(k)\right\}_{k=1}^\infty$
of positive integers
such that, for each $i$,
\begin{gather}\label{zndn}
\sum_k d(x_k)\left(G_{k,i}+1\right)<\infty\quad
\text{with $G_{k,i}=\norm{\ksk{g_k}^{(i)}}_{C(\R,\C)}$}\quad (k\in\Na),\\
\label{sumrk}
\sum_k kz\sk{e_k^\perp}H_{n(k),i}<\infty.
\end{gather}
\end{enumerate}
\begin{rem}
Let $\{ u_n\}$ be an orthonormal basis for $L_2$ such that, for
each $i$,
\begin{gather}\label{1}
\kskob{u_n}^{(i)}\in C(\R,\C)\quad(n\in\Na),\\
\label{2}\norm{\kskob{u_n}^{(i)}}_{C(\R,\C)}\le D_nA_i\quad(n\in\Na),\\
\label{3}
\sum_kD_{n_k}<\infty,
\end{gather}
where $\{D_n\}_{n=1}^\infty$, $\{A_i\}_{i=0}^\infty$ are sequences
of positive numbers, and $\{n_k\}_{k=1}^\infty$ is a subsequence
of $\Na$. Since $d(e_k)\to0$ as $k\to\infty$, the basis $\{ u_n\}$
satisfies Conditions (a) through (c) with $h_k=u_{n_k}$
($k\in\Na$) and
$\{g_k\}_{k=1}^\infty=\{u_n\}\setminus\{h_k\}_{k=1}^\infty$. The
next section of the present paper exemplifies the existence of
such a basis $\{u_n\}$.
\end{rem}

Let us return to the proof of the theorem.
Let $\left\{x_k^\perp\right\}=\left\{e_k^\perp\right\}
\cup(\{e_k\}\setminus\{x_k\})$, and let
$\{f_n\}=\{x_k\}\cup\left\{x_k^\perp\right\}$.
Define a unitary operator $U_\infty:\H\to L_2$ on the basis vectors by setting
\begin{equation}\label{uaction}
U_\infty x_k^\perp=h_k,\quad  U_\infty x_k=g_k\quad\text{for all $k\in \Na$},
\end{equation}
in the harmless assumption that, for each $k\in\Na$,
\begin{equation}\label{ektohnk}
 U_\infty f_k=u_k,\quad  U_\infty e_k^\perp=h_{n(k)},
\end{equation}
where ${n(k)}$ is just that sequence which occurs in Condition (c).

\subsection*{Verifying}
To show that the unitary operator $ U_\infty $ defined in \eqref{uaction} has the
desired properties, it is to be proved that $T= U_\infty S U_\infty ^{-1}$ is an integral
operator having a $SK^\infty$-kernel. For this purpose,
verify that so are the operators $P= U_\infty Q U_\infty ^{-1}$ and $F= U_\infty J^* U_\infty ^{-1}$.

 Using  \eqref{qoper}, \eqref{ektohnk}, one can write
\begin{equation} \label{P}
Pf=\sum_k \ska{f,T^*h_{n(k)}} h_{n(k)}\quad\text{for all $f\in D_T= U_\infty D_S$},
\end{equation}
where, by \eqref{lamgam},
\begin{multline}\label{Thk}
T^*h_{n(k)}=\sum\limits_n\hska{S^*e_k^\perp,f_n} u_n
\\
=
kz\sk{e_k^\perp}\sum\limits_n\hska{e_k^\perp,\varGamma^*f_n} u_n \quad (k\in\Na).
\end{multline}
Prove that, for any fixed $i$, the series
\begin{equation*}\label{imtuko} 
  \sum\limits_n\hska{ e_k^\perp, \varGamma^*f_n} \ksk{u_n}^{(i)}(s) \quad
   (k\in\Na)
\end{equation*}
converge in the norm of $C(\R,\C)$.
Indeed, all these series are pointwise dominated on $\R$
by one series
$$
  \sum\limits_n\hnorm{\varGamma^*f_n}\abs{\ksk{u_n}^{(i)}(s)},
$$
which converges uniformly in $\R$ because its component subseries
\begin{gather*}
  \sum\limits_k\hnorm{\varGamma^*x_k}\abs{\ksk{g_k}^{(i)}(s)},\\
  \sum\limits_k\hnorm{\varGamma^*x_k^\perp}\abs{\ksk{h_k}^{(i)}(s)}
\end{gather*}
are in turn dominated by the convergent series
\begin{equation*}\label{domser} 
 \sum\limits_k d(x_k)G_{k,i}, \quad
 \sum\limits_k\norm{\varGamma^*}H_{k,i},
\end{equation*}
respectively (see \eqref{uaction}, \eqref{defdknums}, \eqref{zndn}, \eqref{hki}).
Whence it follows via \eqref{Thk}
that, for each $k\in\Na$,
\begin{equation}\label{supest}
\norm{\ksk{T^*h_{n(k)}}^{(i)}}_{C(\R,\C)}\le C_ikz\sk{e_k^\perp},\quad
\end{equation}
with a constant $C_i$ independent of $k$.
From \eqref{Thk} it follows also that
\begin{equation}\label{sest}
\norm{T^*h_{n(k)}}\le kz\sk{e_k^\perp}\|\varGamma\|\quad(k\in\Na).
\end{equation}
Consider functions $\P:\R^2\to\C$, $\p:\R\to L_2$, defined, for all
$s$, $t\in\R$, by
\begin{equation} \label{Pp}
\begin{gathered}
\P(s,t)=\sum\limits_k\ksk{h_{n(k)}}(s)
         \overline{\ksk{T^*h_{n(k)}}(t)},\\
\p(s)=\overline{\P(s,\cdot)}=\sum\limits_k
\overline{\ksk{h_{n(k)}}(s)}T^*h_{n(k)}.
\end{gathered}
\end{equation}
The termwise differentiation theorem implies
that, for each $i$ and each integer $j\in[0,+\infty)$,
\begin{gather*}
\dfrac{\partial^{i+j}\P}{\partial s^i\partial t^j}(s,t)
=\sum\limits_k\ksk{h_{n(k)}}^{(i)}(s)
                             \overline{\ksk{T^*h_{n(k)}}^{(j)}(t)},\\
\dfrac{d^i\p}{ds^i}(s)=
\sum\limits_k\overline{\ksk{h_{n(k)}}^{(i)}(s)}T^*h_{n(k)},
\end{gather*}
since, by \eqref{supest}, \eqref{sest}, and \eqref{sumrk},
the series displayed converge (absolutely) in $C(\R^2,\C)$, $C(\R,L_2)$, respectively.
Thus, $\dfrac{\partial^{i+j}\P}{\partial s^i\partial t^j}\in C(\R^2,\C)$,
and $\dfrac{d^i\p}{ds^i}\in C(\R,L_2)$.
Observe also that, by \eqref{sest}, \eqref{sumrk}, and \eqref{Pp}, the series \eqref{P}
(viewed, of course, as one with terms belonging to $C(\R,\C)$) converges
(absolutely) in $C(\R,\C)$-norm  to the function
$$
\kskob{Pf}(s)\equiv\skoba{f,\p(s)}\equiv\int_{\R}\P(s,t)f(t)\,dt.
$$
Thus, $P$ is an integral operator, and $\P$ is its $SK^\infty$-kernel.

Since $\hnorm{S^*e_k}=\hnorm{Je_k}$ for all $k$, from \eqref{null-seq} it follows
via \eqref{schw} that the operator $B$ defined in \eqref{B} is nuclear,
and hence
\begin{equation}\label{snums}
\sum_n s_{n}^{\frac1{2}}<\infty.
\end{equation}
Then, according to \eqref{ES} and \eqref{B}, a kernel which induces
the nuclear operator $F$ can be represented by the series
\begin{equation}\label{fkernels}
\sum_n s_n^{\frac12} U_\infty B^*q_n(s)
\overline{ U_\infty Bp_n(t)}
\end{equation}
convergent almost everywhere in $\R^2$. The functions used
in this bilinear expansion can be written as the series convergent in $L_2$:
$$
   U_\infty Bp_k=\sum\limits_n\hska{ p_k,B^*f_n} u_n, \quad
   U_\infty B^*q_k=\sum\limits_n\hska{q_k,Bf_n} u_n\quad (k\in\Na).
$$
Show that, for any fixed $i$, the functions
$\ksk{ U_\infty Bp_k}^{(i)}$, $\ksk{ U_\infty B^*q_k}^{(i)}$ ($k\in\Na$)
make sense,
are all in $C(\R,\C)$, and their $C(\R,\C)$-norms
are bounded independent of $k$. Indeed, all the series
\begin{equation*} 
  \sum\limits_n\hska{p_k,B^*f_n}\ksk{u_n}^{(i)}(s),\quad
  \sum\limits_n\hska{q_k,Bf_n}\ksk{u_n}^{(i)}(s)\quad (k\in\Na)
\end{equation*}
are dominated by one series
$$
  \sum\limits_n(\hnorm{B^*f_n}+\hnorm{Bf_n})\abs{\ksk{u_n}^{(i)}(s)}.
$$
This series converges uniformly in $\R$, since it consists of two
uniformly convergent in $\R$ subseries
\begin{equation*}
\begin{gathered}
  \sum\limits_k\sk{\hnorm{B^*x_k}+\hnorm{Bx_k}}\abs{\ksk{g_k}^{(i)}(s)},\\
  \sum\limits_k\sk{\hnorm{B^*x_k^\perp}+\hnorm{Bx_k^\perp}}
  \abs{\ksk{h_k}^{(i)}(s)},
\end{gathered}
\end{equation*}
which are dominated by the following convergent series
$$
\sum\limits_k d(x_k)G_{k,i}, \quad
\sum\limits_k 2\|B\|H_{k,i},
$$
respectively (see \eqref{defdknums}, \eqref{schw}, \eqref{zndn}, \eqref{hki}).
Thus, for  functions $F:\R^2\to \C$, $f:\R\to L_2$, defined by
\begin{equation*}
\begin{gathered}
\F(s,t)=\sum_n s_n^{\frac12}\ksk{ U_\infty B^*q_n}(s)\overline{\ksk{ U_\infty Bp_n}(t)},\\
\f(s)=\overline{\F(s,\cdot)}=\sum_n s_n^{\frac12}\overline{\ksk{ U_\infty B^*q_n}(s)} U_\infty Bp_n, 
\end{gathered}
\end{equation*}
one can write, for all non-negative integers $i$, $j$ and all $s$, $t\in\R$,
\begin{equation*}
\begin{gathered}
\dfrac{\partial^{i+j}\F}{\partial s^i\partial t^j}(s,t)
=\sum_n s_n^{\frac12}\ksk{ U_\infty B^*q_n}^{(i)}(s)\overline{\ksk{ U_\infty Bp_n}^{(j)}(t)},\\
\dfrac{d^i\f}{ds^i}(s)
=\sum_n s_n^{\frac12}\overline{\ksk{ U_\infty B^*q_n}^{(i)}(s)} U_\infty Bp_n,
\end{gathered}
\end{equation*}
where  the series converge in $C(\R^2,\C)$, $C(\R,L_2)$, respectively,
because of \eqref{snums}. This implies that $\F$
is a $SK^\infty$-kernel of $F$.

In accordance with \eqref{splitting} and \eqref{ES}, we have,
for each $f\in D_T= U_\infty D_S$,
$$
(Tf)(s)=\int_{\R} \P(s,t)f(t)\,dt+\int_{\R} \F(s,t)f(t)\,dt=
\int_{\R}(\P(s,t)+\F(s,t))f(t)\,dt
$$
for almost every $s$ in $\R$. Therefore $T$ is a
Carleman operator, and that kernel $\K$ of $T$, which is defined by
$\K(s,t)=\P(s,t)+\F(s,t)$ ($s$, $t\in\R$), inherits
the $SK^m$-kernel properties (i), (ii) from its terms, for each $m$.
Consequently, $\K$ is a $SK^\infty$-kernel. The proof of the theorem is complete.

\section{Example of the basis $\{u_n\}$}

What follows is an example of a basis satisfying \eqref{1} through \eqref{3}.
Let $u$ be   Lemari\'e-Meyer wavelet,
$$
u(s)=\dfrac1{2\pi}\int_{\R}e^{i\xi(\frac12+s)}
\sgn\xi b(|\xi|)\,d\xi\quad (s\in\R),
$$
with the bell function $b$ belonging to
$C^\infty(\R)$ (for construction of the Lemari\'e-Meyer wavelets
we refer to \cite{LeMe}, \cite[\S~4]{Ausch}, \cite[Example D, p.~62]{Her}).
Then $u$ belongs to the Schwartz class $\Ss(\R)$, and hence
all the derivatives $\kskob{u}^{(i)}$
are in $C(\R,\C)$.
The ``mother function'' $u$ generates an orthonormal basis for $L_2$ by
$$
u_{jk}(s)=2^{\frac j2}u(2^js-k)\quad  (j,\,k\in\Zi).
$$
 Rearrange, in a completely arbitrary manner, the orthonormal set
$\{u_{jk}\}_{j,\,k\in\Zi}$ into a simple sequence,
so that it becomes $\{ u_n\}_{n\in\Na}$. Since, in view of this rearrangement,
to each $n\in\Na$ there corresponds a unique pair of integers
$j_n$, $k_n$, and
conversely, we can write, for each $i$,
$$
\norm{\kskob{u_n}^{(i)}}_{C(\R,\C)}=\norm{\kskob{u_{j_nk_n}}^{(i)}}_{C(\R,\C)}
\le D_nA_i,
$$
where
$$
D_n=\begin{cases}
2^{j_n^2}&\text{if $j_n>0$,}\\
\left(\dfrac1{\sqrt{2}}\right)^{\abs{j_n}}&\text{if $j_n\le0$,}
\end{cases}
\qquad A_i=2^{\left(i+\frac12\right)^2}\norm{\kskob{u}^{(i)}}_{C(\R,\C)}.
$$
Whence it follows that if $\{n_k\}_{k=1}^\infty\subset\Na$ is a subsequence
such that $j_{n_k}\to -\infty$ as $k\to\infty$, then
$$\sum_kD_{n_k}<\infty.$$
Thus,  the basis $\{u_n\}$ satisfies Conditions \eqref{1} through \eqref{3} and,
consequently, Conditions (a) through (c) in the proof of Theorem.

\section*{Concluding remark}
We conclude this paper by claiming that statements {\bf(A)}, {\bf(B)}, {\bf(C)}
of \propref{msmooth} hold if we replace everywhere the number $m$ by $\infty$.

\section*{Acknowledgements} The author thanks the Mathematical Institute
of the  University of Jena for its hospitality, and specially
Prof. W.~Sickel and Prof. H.-J.~Schmei\ss er for fruitful
discussion on applying wavelets in integral representation theory.

\bibliographystyle{amsplain}
\bibliography{paper}

\end{document}